\documentclass{article}
\usepackage{amssymb,amsmath, amsthm}

\newtheorem{prop}{Proposition}

\newtheorem{lem}{Lemma}
\newtheorem{thm}{Theorem}
\newtheorem{cor}{Corollary}
\newtheorem{rem}{Remark}

\newcommand{\beq}{\begin{equation}}
\newcommand{\eeq}{\end{equation}}

\def\eq#1/{(\ref{e:#1})}

\makeatletter
\@addtoreset{equation}{section}
\makeatother
\makeatletter
\@addtoreset{equation}{section}
\makeatother

\newcommand{\ta}{\theta}

\newcommand{\si}{\sigma}

\newcommand{\ff}{\infty}
\newcommand{\ra}{\rightarrow}

\newcommand \QED{\hbox{\hskip6pt\vrule height 8pt width 6pt}}

\newcommand {\vxn} {{\bf{X}}_{n,r(n)}}

\def\Real{\mathrm{I\!R}}

\begin{document}
\title{Weak Convergence Results for Multiple Generations of a Branching Process}
\author{James Kuelbs \\ Department of Mathematics\\   University of Wisconsin  \\ Madison, WI 53706-1388 \\ E-mail: kuelbs@math.wisc.edu \and Anand N. Vidyashankar$^{\star}$ \\Department of Statistical Science \\Cornell University \\Ithaca, NY  14853-4201\\ E-mail: anv4@cornell.edu} 
\date{ }
\maketitle
\begin{abstract}
We establish limit theorems involving weak convergence of multiple generations
of critical and supercritical branching processes. These results arise naturally when dealing with 
 the joint asymptotic behavior of 
functionals defined in terms of several generations of such processes. Applications of  our main result include a functional central limit theorem (CLT), a Darling-Erd\"os result, and an extremal process result.
The limiting process for our functional CLT is an infinite dimensional Brownian motion with sample paths in the infinite product space $(C_0[0,1])^{\infty}$, with the product topology, or in Banach subspaces  of $(C_0[0,1])^{\infty}$  determined by norms related to the distribution of the population size of the branching process. As an application of this CLT we obtain a central limit theorem for ratios of weighted sums of generations of a branching processes, and also to various maximums of these generations. The Darling-Erd\"os result and the application to extremal distributions also include infinite dimensional  limit laws. Some branching process examples where the CLT fails are also included. 

\hfill\\\rule{\textwidth}{0.1pt}
{$^{\star}$ Research Supported in part by a grant from NSF DMS 000-03-07057 and also by grants from the NDCHealth Corporation}
{\em AMS}~1991~ {\em Subject Classification:}~60J80 60F17 60B10 60B12 \\
{\em Short title:} {Multiple Generations of Branching Processes}
\end{abstract}
\vfill\eject
\section{\bf Introduction}
Our interest in limit theorems for multiple generations of branching processes is motivated by both practical and theoretical considerations. The practical side stems from the use of branching processes to model certain aspects of scientific experiments. One such problem area is  Polymerase Chain Reaction (PCR) experiments.  In such an experiment, an initial amount of DNA is amplified for use in various biological experiments. The PCR experiment evolves in three phases;  an exponential phase, a linear phase, and a pleateau phase, with  branching processes and their variants frequently used to model 
the exponential phase.
 One of the goals of such experiments is to ``quantitate'' the initial number of DNA molecules in a sample or equivalently, estimate the number of ancestors in a branching process \cite{Jacob1}. The statistical estimate of the initial number of ancestors is a function of the estimate of the mean of the branching process \cite{Jacob1}, and in order to make this estimate, data are used from the last few cycles (generations) at 
the end of the exponential phase.  Since the cycle(generation) corresponding to the end of exponential phase is somewhat arbitrary, it is natural to consider the joint distributions of the generations involved to determine
whether two different scientists with different choices for the end of the exponential phase obtain consistent results. Furthermore, these joint distributions can also be used to estimate the end of the exponential phase.

Theoretical motivation for our results involves the desire to understand analogues of classical functional limit theorems for i.i.d. sequences that hold for multiple generations of the stochastic processes arising in the branching setting. What we present here deals with weak convergence results. Theorem 1 is our main result, and allows a large number of applications, a few of which are presented explicitly as Applications 1-3, and Theorems 2 and 3 in Section 2. Application 1 is a functional CLT, yielding a Donsker type result, Application 2 a Darling-Erd\"os result, and Application 3 an extremal process result, all obtained under best possible conditions. For example, in the functional CLT we use only second moments, and in the Darling-Erd\"os result we use the moment condition shown in \cite{EU} to be necessary for this result for i.i.d. sequences. A similar comment applies to the application to extremal processes. Here the regularly varying tail condition assumed is precis
 ely that required for the limiting maximal distribution at $t=1$ to exist for an i.i.d sequence. Other applications are also possible once one has Theorem 1 available, but in Theorems 2 and 3 we turn to some applications of our functional CLT. Theorem 2 yields a strengthening of the functional CLT to the Banach spaces
$c_{0,\lambda}(C_0[0,1])$. 
Another consequence of Application 1 is a new proof of the CLT for the non-parametric maximum likelihood estimate of the mean of a supercritical branching process. A previous proof of this in \cite{S} involves a martingale CLT, whereas the proof herein is an elegant application of our functional CLT result with $t=1,$ and the asymptotic independence obtained in the coordinates of the limiting process. Moreover, our proof allows us to extend this result to allow the application of a broad range of weights on the various generations. In \cite{S} all the weights are equal to one. 

In order to describe our results in more detail we begin with a brief description of the branching process.  Let $ \{\xi_{n, j}, j \ge 1, n \ge 1\}$ denote a double array of integer valued i.i.d. random variables defined on the probability space $(\Omega, \mathcal{F}, P)$, and having probability distribution $\{p_j: j \ge 0\}$, i.e.
$P(\xi_{1,1}=k)= p_k. $
Then $\{ Z_n: n \ge 0 \}$ denotes the Galton-Watson process initiated by a single ancestor $Z_0 \equiv 1$. It is iteratively defined on $(\Omega, \mathcal{ F}, P)$ for  $n \ge 1$ by
$$
Z_{n}= \sum_{j=1}^{Z_{n-1}}\xi_{n, j}. 
$$
Let $m = E(Z_1)$.  It is well known that if $m >1$ ({\it i.e.} the process is supercritical), then $Z_n \ra \ff$ with positive probability and that the probability that the process becomes extinct, namely $ q $, is less than one.  The complement of the set  $\cup_{n=1}^{\infty} \{Z_n =0\}$ is the so called survival set, and is denoted by $S$. If $m >1$, then $P(S)=1-q$ and $Z_n \rightarrow \infty$  a.s. on $S$. 
Also, $q =0$ if and only if $p_0 =0$.  If $m \le 1$, then assuming $p_1 \neq 1$ when $m=1$, the process becomes extinct with probability one, i.e. $P(S)=0$. To avoid degenerate situations we assume throughout the paper that $p_0+p_1<1.$

The paper is organized as follows: Section 2 develops the basic notation and states the main results of the paper.  Section 3 contains the proof of Theorem 1, and Sections 4 and 5 that of the CLT applications in Theorem's 2 and 3, respectively. Section 6 contains examples providing some insight into the CLT for subcritical processes, and also for supercritical processes when one uses deterministic normalizations. In this latter example one does not get a Gaussian limit law, but a certain mixture of Gaussian laws. This mixture can be anticipated from the Kesten-Stigum result, but its precise expression requires some interesting analysis. In particular, these examples show precisely why the random normalizations used in our theorems are possibly the "best choice" if one wants classical results to persist in limit theorems for multiple generations of these processes.   

{\bf Acknowledgment.} The authors thank the referee for a very careful reading of the manuscript, and for making several important suggestions that led to improvements in the paper.
In particular, Lemma 3 as presented here is an elegant modification of our original arguments, and is based on the referee's comments.

\section{Notation and Main Results} 
In this section we state the main result of the paper. This result allows us to obtain a wide variety of limit theorems for branching processes based on $r(n)$-generations, where $1 \leq r(n) \leq n$. Following its statement we present some interesting consequences and applications. In particular, in these applications the integer sequence $\{r(n)\}$ may approach infinity as n goes to infinity. As will be seen, they all follow rather immediately from our main result when combined with various classical limit theorems for i.i.d. sequences.

Throughout $(M,d)$ is a complete separable metric space with distance $d$, and  $M^{\infty}$ denotes the infinite product of copies of $M$ with the product topology, metrized by
\begin{align}{\label{def:metricforPT}}
d_{\ff}({\bf x}, {\bf y}) = \sum_{k \ge 1}\frac{1}{2^k}\frac{d(x_k,y_k)}{1+d(x_k,y_k)}.
\end{align}
In our applications $M$ is the real line or some function space. If $M$ is the real line, then the distance is the usual one, and for our functional CLT application $M$  denotes the set of all continuous functions on [0,1] that vanish at 0, which we denote by $C_0[0,1]$. Of course, then  $C_0[0,1]$ is a Banach space in the supremum norm 
\begin{align}
q(f)= \sup_{0 \leq t \leq1}|f(t)|,
\end{align}
and the distance used is $d(f,g)=q(f-g), f,g \in C_0[0,1]$. Application 3 below contains a different choice of M, and others are certainly possible, but these suffice to provide a sampling of possible consequences of our main theorem.

Since we want to study the asymptotic behavior of $r(n)$ generations of the branching process, and  $r(n)$ may well converge to infinity, it is useful for these purposes to define

\begin{align}{\label{def:vector-process}}
\vxn  \equiv (X_{n, Z_{n-1}}, X_{n-1, Z_{n-2}}, \cdots X_{n-r(n)+1, Z_{n-r(n)}}, z, z, \cdots ),
\end{align}
where $z$ is a fixed element in $M$,
\begin{align}
X_{n-j+1, Z_{n-j}}= H_{Z_{n-j}}(\xi_{n-j+1,1}, \cdots,\xi_{n-j+1,Z_{n-j}}),
\end{align}
and
the mappings $H_k(\cdot), k \geq 1,$ take $R^k$ into $M$ are Borel measurable. We also define $X_{n-j+1, 0}=z$ for $1 \le j \le r(n) \le n$.  Since in our results we condition on $Z_{n-1} >0$, the choice of $z$ is immaterial. Hence  $\vxn$ is an element of the infinite product space $M^{\ff}$. Moreover,
in our applications $M$ always contains a zero element which we denote by $0$, and if we take the fixed element $z \in M$ in (2.3) to be this $0$, then we have
\begin{align}{\label{def:vector-process}}
\vxn  \equiv (X_{n, Z_{n-1}}, X_{n-1, Z_{n-2}}, \cdots X_{n-r(n)+1, Z_{n-r(n)}}, 0, 0, \cdots ).
\end{align}

We will use $\Rightarrow$ to denote weak convergence of probability measures. Our main theorem for the random vectors $\vxn$ is the following. 

\begin{thm}{\label{thm:weakconvergence}} Let $m \geq 1>p_1,$  assume $1 \leq r(n) \leq n$ with  $\lim_{n \ra \ff} r(n)= \ff$, and that $\vxn$ is defined as in (2.3)-(2.4). Also assume that  if $\{\xi_j: j \geq 1\}$ are i.i.d. non-negative integer valued random variable with $\mathcal{L}(\xi_1)=\mathcal{L}(Z_1)$, then the $M$-valued random elements $\{H_k: k \geq 1\}$ used to define $\vxn$ are such that
\begin{align}
H_k(\xi_1,\cdots,\xi_k) \Rightarrow H 
\end{align}
on $(M,d)$.
Moreover, assume that 
\begin{align}
\lim_{n \rightarrow \infty} \frac{P(Z_{n-1} >0)}{P(Z_n >0)} = 1. 
\end{align}
Then the  probability measures 
\begin{align}
\mu_n= \mathcal{L}(\vxn | Z_{n-1}>0)
\end{align}
converge weakly on $(M^{\ff},d_{\ff})$, i.e. we have
\begin{align}
\mu_n \Rightarrow \mathcal{L}(B_1,B_2,\cdots),
\end{align}
where the $B_i$'s  are independent copies of $H$.
\end{thm}

\bigskip

\begin{rem} 
If $m>1$, then $\lim_{n \rightarrow \infty} P(Z_n>0)= 1-q>0$ and hence the condition (2.7) holds. If $m=1$ and $0<\sigma^2 =E((Z_1 -m)^2) < \infty$, then (2.7) also follows from Theorem 1 in [1, p19]. Furthermore, as was pointed out by the referee, by Lemma 2 in  \cite{SL}, there are interesting examples where $m=1$ and $\sigma^2=\infty$, yet (2.7) holds.
\end{rem}

\bigskip

Next we present three immediate applications of Theorem 1. They include a functional CLT, a Darling-Erd\"os Theorem, and also an extremal process result. It is interesting to observe that the limiting distributions of the coordinates of $\vxn$ are asymptotically independent, whereas the generations of the branching process itself are correlated. 
\bigskip

{\bf Application 1}: Let $m \geq 1$, $0< \si^2=E((Z_1-m)^{2}) < \ff$,  and assume $1 \leq r(n) \leq n$ with  $\lim_{n \ra \ff} r(n)= \ff$. Take $M=C_0[0,1]$ with the sup-norm $q$, define  $H_k(\xi_1,\cdots,\xi_k)(0)=0$, and for $0 \leq t \leq 1$ set
\begin{align}
H_k(\xi_1,\cdots,\xi_k)(t)= \frac{1}{\sigma \sqrt k} \sum_{i=1}^{\lfloor{tk}\rfloor}(\xi_i-m) + (tk -\lfloor{tk}\rfloor)\frac{1}{\si\sqrt k}(\xi_{\lfloor{tk}\rfloor+1}-m).
\end{align}
Then Donsker's Invariance principle implies (2.6) holds with $\mathcal{L}(H)$  the probability measure induced on $C_0[0,1]$ by a standard Brownian motion starting at zero when $t=0$. If $\vxn$ is defined as  in (2.3-5) with $H_k$ as in (2.10), and $(C_0[0,1])^{\ff}$ has the product topology induced when using the norm $q$ on $C_0[0,1]$, then an immediate consequence of Theorem 1 is that the  probability measures 
\begin{align}
\mu_n= \mathcal{L}(\vxn | Z_{n-1}>0)
\end{align}
converge weakly there, i.e. we have
\begin{align}
\mu_n \Rightarrow \mathcal{L}(B_1,B_2,\cdots),
\end{align}
where the $B_i$'s  are independent Brownian motions.
\bigskip
\bigskip

{\bf Application 2}: Let $m \geq 1$, $0< \si^2=E((Z_1-m)^{2}) < \ff$, $\lim_{t \rightarrow \ff} LLtE(Z_1^2I(|Z_1| \geq t)) =0$,  where $Lt= \log_e(t\vee e)$ and $LLt=L(Lt)$.  In addition, assume $1 \leq r(n) \leq n$ with  $\lim_{n \ra \ff} r(n)= \ff$, and take $M=R^1$. Define 
 \begin{align}
H_k(\xi_1,\cdots,\xi_k)= a_k \max_{1 \leq j \leq k}\frac{ \sum_{i=1}^{j}(\xi_i-m) }{\sigma \sqrt j}-b_k,
\end{align}
where $a_k =(2LLk)^{\frac{1}{2}}$ and $b_k= 2LLk  + \frac{1}{2}LLLk - \frac{1}{2}L(4\pi)$.
Then the Darling-Erd\"os Theorem as in Theorem 2 of \cite{EU} implies (2.6) holds with $\mathcal{L}(H)$  the probability measure induced on $M$ by the cumulative distribution function
 \begin{align}
G(x) = \exp\{-e^{-x}\}, ~ -\infty<x<\ff.
\end{align}
If $\vxn$ is defined as  in (2.3-5) with $H_k$ as in (2.13), and $M^{\ff}=R^{\ff}$ has the product topology , then an immediate consequence of Theorem 1 is that the  probability measures 
\begin{align}
\mu_n= \mathcal{L}(\vxn | Z_{n-1}>0)
\end{align}
converge weakly there, and we have
\begin{align}
\mu_n \Rightarrow \mathcal{L}(B_1,B_2,\cdots),
\end{align}
where the $B_i$'s  are independent random variables with cumulative distribution function $G(x)$ given by (2.14).
\bigskip
\bigskip

{\bf Application 3}: Let $m \geq 1$,  assume $F(x) = P(Z_1\leq x)<1$ for all $ x\in R^1$, and that $1-F(x)$ is regularly varying at $\ff$ with exponent $-\alpha$ where $\alpha> 1.$  In addition, assume $1 \leq r(n) \leq n$ with  $\lim_{n \ra \ff} r(n)= \ff$, and that (2.7) holds. Then by
Theorem 6.3, p.455 of  \cite{G}, there exists $a_j>0$ such  that
\begin{align}
\lim_{j \rightarrow \ff}P( \frac{1}{a_j} \max\{0,\xi_1,\cdots,\xi_j\} \leq x) = \exp\{- x^{-\alpha}\}, ~x>0, 
\end{align}
and zero for $x \leq 0$.
Now define the extremal process $m_k(t)$ which is zero in $[0,\frac{1}{n})$ and
\begin{align}
m_k(t) = \frac{1}{a_k} \max\{0,\xi_1,\cdots,\xi_j\}, ~\frac{ j}{n} \leq t <\frac{ j+1}{k},~ j=1,\cdots,k-1,
\end{align}
and 
\begin{align}
m_k(t) =  \frac{1}{a_k} \max\{0,\xi_1,\cdots,\xi_k\},~ t\geq 1.
\end{align}
Let $M$ denote the finite, non-decreasing functions $z(t)$on $[0,\ff)$ such that $z(0)=0$ and $ z(t)=z(1)$ for  $t \geq 1$.Then $M$ is a complete separable metric space in the Levy metric $d_L$ on $M$, and if
\begin{align}
H_k(\xi_1,\cdots,\xi_k)(t)= m_k(t),~ 0\leq t <\ff,
\end{align}
then by Theorem 2.1 and 3.1 of \cite{L} we have (2.6) where $\{H(t):0\leq t <\ff\}$ is a Markov extremal process with sample paths in $M$. Therefore,
if $\vxn$ is defined as  in (2.3-5) with $H_k$ as in (2.20), and $M^{\ff}$ has the product topology, then an immediate consequence of Theorem 1 is that the  probability measures 
\begin{align}
\mu_n= \mathcal{L}(\vxn | Z_{n-1}>0)
\end{align}
converge weakly there, and we have
\begin{align}
\mu_n \Rightarrow \mathcal{L}(B_1,B_2,\cdots),
\end{align}
where the $B_i$'s  are independent copies of the Markov process $\{H(t):0\leq t <\ff\}$.
\bigskip

As is easily seen, Theorem 1 combined with other classical limit theorems for i.i.d sequences provides many possible limit theorems for suitable choices of the random elements $\vxn$. However, what we turn to next are some applications and extensions of the functional CLT of application one.
The first involves a functional CLT in  Banach subspaces of $(C_0[0,1])^{\ff}$ determined by weighted analogues of the $q$-norm. That is, let $\lambda =\{\lambda_j : j \geq 1\}$ be a sequence of strictly positive numbers, and for ${\bf f}=(f_1,f_2,\cdots) \in (C_0[0,1])^{\infty}$ define
\begin{align}
q_{\lambda}({\bf f}) = \sup_{j \geq 1} \lambda_j||f_j||,
\end{align}
where $||\cdot||$ is the supremum norm on $C_0[0,1]$. Also, let $c_{0,\lambda}(C_0[0,1])$ be the subspace of $(C_0[0,1])^{\infty}$ given by
\begin{align}
c_{0,\lambda}(C_0[0,1]) = \{{\bf f}=(f_1,f_2,\cdots) \in (C_0[0,1])^{\infty}: \lim_{j \rightarrow \infty} \lambda_j||f_j||=0\}.
\end{align}
Then $q_{\lambda}({\bf f})$ is a norm making the subspace $c_{0,\lambda}(C_0[0,1])$ a Banach space.
 
 \bigskip

As before we will use $\Rightarrow$ to denote weak convergence of probability measures. Our functional central limit theorem in $c_{0,\lambda}(C_0[0,1])$ is the following. 

\begin{thm}{\label{thm:weakconvergence}} Let $m \geq 1$ and assume $1 \leq r(n) \leq n$ with  $\lim_{n \ra \ff} r(n)= \ff$. Also assume that  the offspring distribution $\mathcal{L}(\xi)=\mathcal{L}(Z_1)$ is such that $0< \si^2=E((\xi-m)^{2}) < \ff$ and satisfies one of the following conditions:
\begin{align}
P(|\xi-m|\ge x) \le \beta e^{-\theta x^2},~ {\rm{for~ all}
} ~ x \ge 0,~\rm {or}
\end{align}
\begin{align}
E(|\xi-m|^{\rho}) < \ff~\rm{for~ some} ~\rho\ge 2, 
\end{align}
and that $r(n) = o(n)$. Let $\vxn$ be defined as in (2.3-5)
with $H_k(\xi_1,\cdots,\xi_k)$ as in (2.10).
If (2.25) holds and we take $\lambda = \{\lambda_j\}$ where $\lambda_j= (\delta_j\log(j+3))^{-\frac{1}{2}}$ and $\lim_{j \rightarrow \infty} \delta_j= \infty$, then 
on the Banach space $c_{0,\lambda}(C_0[0,1])$ the probability measures 
\begin{align}
\mu_n= \mathcal{L}(\vxn | Z_{n-1}>0).
\end{align}
are such that 
\begin{align}
\mu_n \Rightarrow \mathcal{L}(B_1,B_2,\cdots),
\end{align}
where the $B_i$'s  are independent standard Brownian motions.
If instead we assume (2.26) and $\lambda = \{\lambda_j\}$ where $\lambda_j= j^{-\frac{(1+\delta)}{\rho}}$ and $\delta>0$, then we again have (2.28) on $c_{0,\lambda}(C_0[0,1])$. 
 \end{thm}

\begin{rem}
The condition  $r(n) = o(n)$ in Theorem 2 can be weaken somewhat. For example, if $m=1$ and  $0< \si^2=E((\xi-m)^{2}) < \ff$ we need only assume that $\limsup_{n \rightarrow \infty} \frac{r(n)}{n}<1$. Similarly, if $m>1$, we may replace $r(n)= o(n)$ by $\lim_{n \rightarrow \infty} (n - r(n)) = \infty$. Both of these improvements follow from refinements of (4.7)-(4.8) of Lemma 4 below, and were pointed out by the referee. We emphasize that throughout the assumption $\lim_{n \rightarrow \infty} r(n) = \infty$ is in effect.
\end{rem}

If $G$ is Gaussian random variable with mean zero and variance one, then for all $ x \geq 0$
$$
P(\sup_{0 \leq t \leq 1}B(t) \leq x) = P(|G| \leq x).
$$
Hence Theorem 2 and the continuous mapping theorem applied to the processes $\{\vxn(\cdot): n \geq 1\}$ with values in $c_{0,\lambda}(C_0[0,1])$ immediately imply the following result. 
\bigskip

\begin{cor} If (2.25) or (2.26) holds with corresponding $\{ \lambda_j: j \geq 1\}$ as indicated, then the conditions of Theorem 2 imply that
$$
\lim_{n \rightarrow \infty} P(\max_{ 1 \leq j \leq r(n)} \lambda_j\frac{ (Z_{n-j+1} - mZ_{n-j})}{\sigma Z_{n-j}^{\frac{1}{2}}} \leq x|Z_{n-1}>0) = P(\sup_{j \geq 1} \lambda_jG_j \leq x),
$$
where $G_1,G_2,\cdots$ are i.i.d. $N(0,1)$ random variables.
In addition, we also have
$$
\lim_{n \rightarrow \infty} P(\max_{ 1 \leq j \leq r(n)} \lambda_j( \frac{\max_{1 \leq k \leq Z_{n-j}} (0 \vee \sum_{i=1}^k(\xi_{n-j+1,i}-m))}{\sigma Z_{n-j}^{\frac{1}{2}}}) \leq x|Z_{n-1}>0)= P(\sup_{j \geq 1} \lambda_j|G_j| \leq x).
$$
\end{cor}
\bigskip

\begin{rem} If one wants results similar to those of the corollary with $\lambda_j=1$, then Theorem 2, or application one applies, as long as we restrict the maximum to be over only finitely many  $j$'s, say $ j \in \{1,2,\cdots,d\}$.
\end{rem}
\bigskip


In Theorem 3 below we obtain  a CLT for ratios of weighted 
sums of a supercritical branching process $\{Z_n: n \geq 1\}$. 
When the weights are all one the result appeared in \cite{S} using a martingale CLT for the proof. Our proof is completely different. It uses Application 1 in an important way and allows the ratios to consist of weighted sums. We begin with some notation.
\bigskip

In Theorem 3 we assume $\{b_j: j \geq 1\}$ is a sequence of non-negative numbers 
and set
$$
X_n= (\frac{N_{n}}{D_{n}} -m)\sqrt{D_{n}},
$$ 
where
$$
N_{n}= b_1Z_n + b_2Z_{n-1}+\cdots +b_nZ_{1},
$$
$$
D_{n}= b_1Z_{n-1} + b_2Z_{n-2}+\cdots +b_nZ_{0},
$$
and we understand $X_n$ to be zero if $D_n=0$.
Then we have the following CLT.
 
 \begin{thm} Let $m>1,~0< \si^2 \equiv E((Z_1-m)^2) <\infty$, and assume $\{b_j: j \geq 1\}$ is a sequence of non-negative numbers with 
$$
0< \kappa \equiv \sum_{j=1}^{\infty}\frac{ b_j^2}{ m^{j}} < \infty.
$$
For $k \geq 1$, let
\begin{align}
\theta_k = \frac{b_k}{m^{k}}(\kappa)^{-1},
\end{align}
and set
$$
\Lambda^2= \frac{\sigma^2}{\kappa} \sum_{j=1}^{\infty} \frac{b_j^2}{m^j}.
$$ 
Then,  $\Lambda^2<\infty$, 
 and for all real $x$ we have
\begin{align}
\lim_{n \rightarrow \infty} P(X_n \leq x|S) = P(G \leq x),
\end{align}
where $G$ is a mean zero Gaussian random variable with $E(G^2) = \Lambda^2$. In particular, if $b_j=0$ or $b_j=1$ for all $j\ge 1$ and some $b_j>0$, then $\Lambda^2=\si^2$ and (2.30) holds with $G$  a mean zero Gaussian random variable with $E(G^2) = \si^2$.
\end{thm}



\begin{rem} If we condition on $\{Z_{n-1} >0\}$ instead of $S$ in Theorem 3, the limit is the same. 
\end{rem}


\section{Proof of Theorem 1}

The proof of Theorem 1 is based on a lemma for weak convergence in infinite product spaces, and an iterative technique developed in Lemma 3 below. This iterative lemma also is applicable to the proof of Proposition 1, which appears in Section 6.

Let $(M, d)$ be a complete separable metric space, $\mu$  a Borel probability measure on $(M, d)$, and $\pi: M \ra M$  Borel measurable. Define,
$$
\mu^{\pi}(A)= \mu(\pi^{-1}(A))
$$
for all Borel sets $A$ of $(M, d)$. Let $M^{\ff}$ denote the infinite product space with the product topology and typical point ${\bf s}= (s_1, s_2, \cdots )$. 
If $z$ is a fixed point in $M$, we define the mapping
$\pi_l: M^{\ff} \ra M^{\ff} $, for $l \ge 1$, by
$$
\pi_l({\bf s}) = (s_1, s_2,  \cdots, s_l, z,z, \cdots, ).
$$
We now indicate a lemma concerning weak convergence in product spaces.  Its proof is easily anticipated.
\bigskip

{\begin{lem} Let $M$ be as above and assume $\{\mu_n: n \ge 1\}$ and $\mu_{\ff}$ are Borel probability measures on $M^{\ff}$ with the product topology. Then $\{\mu_n: n \ge 1\}$ converges weakly to $\mu_{\ff}$ if and only if $\mu_n^{\pi_l}$ converges weakly to $\mu_{\ff}^{\pi_l}$ for all $l \ge 1$.
\end{lem}}

\bigskip

The next lemma is used in the proof of Theorem 1. When $m=1$ the result follows from Theorem 2, p. 20, of \cite{AN}, and when $m>1$, a bit of calculation shows that it follows from Lemmas 1 and 2, p. 4-5, of \cite{AN}. 
\bigskip

\begin{lem} Let $\{Z_n: n \ge0\}$ be a Galton-Watson process with $Z_0=1$. If $m\ge 1>p_1$, then  for each $J \in [1,\infty)$
\begin{align}
\lim_{n \rightarrow \infty} \frac{P(1 \leq Z_n \leq J)}{P(Z_n >0)} =0.
\end{align}
\end{lem}

{\bf Proof of (2.9)} Let $\mu$ denote the probability measure induced by $H$ on $M$, and $\mu_{\ff}$ be the infinite product measure formed by $\mu$ on $M^{\ff}$. Also let $\mu_n$ denote the law of $\vxn$ when $Z_{n-1}$ is conditioned to be stricty positive,i.e.
for $A$ a Borel subset of $M^{\ff}$ we have
$$
\mu_n(A)= P(\vxn \in A|Z_{n-1}>0).
$$
By Lemma 1 
it is sufficient to establish, for each $l \ge 1$, the weak convergence of $\mu_{n}^{\pi_l}$  to $\mu_{\ff}^{\pi_l}$. If we identify the range space of $\pi_l$ with $M^l$ in the obvious way, then it suffices to show that on $M^l$ we have that
$$
\lambda_n= \mathcal{L}( X_{n, Z_{n-1}}, X_{{n-1}, Z_{n-2}}, \cdots, X_{{n-l+1}, Z_{n-l}}|Z_{n-1}>0)
$$
converges weakly to $(\mu)^l$, the l-fold product of $\mu$ on that space.

To establish weak convergence of $\lambda_n$ to $( \mu)^l$, it is sufficient by Theorem 2.2 of \cite{BI} to show for arbitrary continuity sets $E_i$ of the measure $\mu$ on $M$ that
\begin{align}
\lim_{n \ra \ff} \lambda_n(E_1 \times E_2 \times \cdots \times E_l) = \prod_{j=1}^{l}\mu(E_j).
\end{align}
We will now verify

{\begin{lem} 
Let $\{Z_n: n \ge0\}$ be a Galton-Watson process with $Z_0=1$, $m\ge 1>p_1$, and (2.7) holding. Also let $\mathcal{F}_0 = \{\phi, \Omega \}$ and $\mathcal{F}_n= \sigma ( \{\xi_{k, j}: j \ge 1\}: 1 \le k \le n)$ for $ n \ge 1$. 
Let $l \geq 1$ be an integer and let $\{ \beta_{n,i}, \gamma_{n,i}: 1\le i \le l <n<\infty\}$ be random variables such that $\beta_{n,i}$ is $\mathcal{F}_{n-l}$ measurable, $\gamma_{n,i}$ is $\mathcal{F}_{n-i+1}$ measurable for all $1 \le i \le l<n$, and for some constant $c>0$ we have $|\beta_{n,i}|\leq c$ a.s. and  $|\gamma_{n,i}|\leq c$ a.s. for all  $1 \le i \le l<n$. Suppose for every $i=1,\cdots,l$ that
\begin{align}
\frac{1}{q_{n-i}}I_{Z_{n-i}>0} E(\gamma_{n,i}- \beta_{n,i}|\mathcal{F}_{n-i}) \rightarrow 0~ \rm {in}~ L^1(P),  
\end{align}
where $q_n=P(Z_n>0)$ for all $n \geq 1$, and set 
\begin{align}
\Gamma_{n,i} = I_{\{Z_{n-i}>0\}}\prod_{\nu=i}^l \gamma_{n,\nu}, ~~B_{n,i}= I_{\{Z_{n-i}>0\}}\prod_{\nu=i}^l \beta_{n,\nu},
\end{align} 
for every $1 \le i \le l<n$. Then for every $j=1,\cdots,l$ we have
\begin{align}
\frac{1}{q_{n-j}}E(\Gamma_{n,j} -B_{n,j}|\mathcal{F}_{n-l}) \rightarrow 0 \rm ~{in}~L^1(P)
\end{align}
and
\begin{align}
\lim_{n \rightarrow \infty} E(\prod_{i=j}^l \gamma_{n,i}- \prod_{i=j}^l \gamma_{n,i}~|Z_{n-j}>0)=0.
\end{align}
If  $\zeta$ maps the non-negative integers into $[0,\ff)$ with $\lim_{k \rightarrow \ff} \zeta(k) =0$, then we also have
 \begin{align}
\lim_{n \rightarrow \infty} E(\zeta(Z_n)|Z_n>0)=0.
\end{align}
\end {lem}}

{\bf Proof.}  First we observe that (3.6) follows from (3.5) by taking expectations.  The proof of (3.5) goes by backwards induction on $j$. By (3.3) with $i=l$, we see that (3.5) holds with if $j=l$, so now let $1 \le j <l$ be a given integer such that  (3.5) holds for $j+1$.

Set $F_k= \{Z_k >0\}$. Then the $F_k$'s decrease as $k$ increases, and
$$
\Gamma_{n,j} - B_{n,j} = I_{F_{n-j}}\gamma_{n,j}\Gamma_{n,j+1} - I_{F_{n-j}}\beta_{n,j}B_{n,j+1},
$$
Hence, if 
$$
\hat \gamma_{n,j}= I_{F_{n-j}} \gamma_{n,j} \rm~ and~\hat \beta_{n,j}= I_{F_{n-j-1}} \beta_{n,j},
$$
we then have
\begin{align}
E((\Gamma_{n,j} - B_{n,j})| \mathcal{F}_{n-l})  = I_n + II_n +III_n,
\end{align}
where
\begin{align}
I_n=E((\hat \gamma_{n,j} - \hat \beta_{n,j})\Gamma_{n,j+1}|\mathcal{F}_{n-l}),
\end{align}
\begin{align}
II_n=E( \hat \beta_{n,j}( \Gamma_{n,j+1}-B_{n,j+1})| \mathcal{F}_{n-l}),
\end{align} 
and
\begin{align}
III_n=E( \hat \beta_{n,j}( B_{n,j+1}-I_{F_{n-j}}B_{n,j+1})| \mathcal{F}_{n-l}).
\end{align}
Therefore, it suffices to show that the three quantities $\frac{I_n}{q_{n-j}}, \frac{II_n}{q_{n-j}}$, and $\frac{III_n}{q_{n-j}}$ converge to zero in $L^1(P)$.

Let
$$
\Delta_{n,j}=E(\hat \gamma_{n,j}-\hat \beta_{n,j} |\mathcal{F}_{n-j}).
$$
Then, since $|\gamma_{n,i}| \leq c$, $\Gamma_{n,j+1}$ is $\mathcal{F}_{n-j}$ measurable, and $1\leq j<l$, we have
\begin{align}
E(|I_n|)= E(|E(E[\Delta_{n,j}|\mathcal{F}_{n-j}]\Gamma_{n,j+1}|\mathcal{F}_{n-l})|)\le \max\{1,c^l\}E(|E[(\hat \gamma_{n,j} - \hat \beta_{n,j})|\mathcal{F}_{n-j}]|)
\end{align}
Since $|\beta_{n,j}| \leq c$, we have $|I_{F_{n-j}} \beta_{n,j} - \hat \beta_{n,j}| \le cI_{F_{n-j-1} \cap F_{n-j}^c},$
and hence (2.7) implies that
\begin{align}
\lim_{n \rightarrow \infty} q_{n-j}^{-1}E(|I_{F_{n-j}} \beta_{n,j} - \hat \beta_{n,j}|)=0.
\end{align}
Hence by combining (3.3) with $i=j$, and (3.13), we see that
\begin{align}
\lim_{n \rightarrow \infty} q_{n-j}^{-1}E(|\Delta_{n,j}|)=0.
\end{align}
Thus $\frac{I_n}{q_{n-j}}$ converges to zero in $L^1(P)$.

Next we observe that
\begin{align}
\beta_{n,j}\Gamma_{n,j+1}=\hat \beta_{n,j}\Gamma_{n,j+1} \rm~{and}~\beta_{n,j}B_{n,j+1}=\hat \beta_{n,j}B_{n,j+1}.
\end{align}
Thus $\beta_{n,j}$ being $\mathcal{F}_{n-l}$ measurable implies
\begin{align}
\frac{II_n}{q_{n-j}} =\beta_{n,j}E((\Gamma_{n,j+1} - B_{n,j+1})|\mathcal{F}_{n-l})/q_{n-j}
\end{align}
and 
\begin{align}
\frac{III_n}{q_{n-j}} =\beta_{n,j}E((B_{n,j+1} - I_{F_{n-j}}B_{n,j+1})|\mathcal{F}_{n-l})/q_{n-j}.
\end{align}
Therefore,
\begin{align}
\frac{E(|II_n|)}{q_{n-j}} \le cE(|E(\Gamma_{n,j+1} - B_{n,j+1})|\mathcal{F}_{n-l})|)/q_{n-j},
\end{align}
and since the induction hypothesis provides (3.5) for $j+1$ and (2.7) holds, we have that (3.18) implies
\begin{align}
\lim_{n \rightarrow \infty} \frac{E(|II_n|)}{q_{n-j}}=0.
\end{align}

Using (3.17), a  similar argument implies that
\begin{align}
\lim_{n \rightarrow \infty} \frac{E(|III_n|)}{q_{n-j}}=0
\end{align}
provided we show
\begin{align}
\lim_{n \rightarrow \infty} \frac{E(|B_{n,j+1} - I_{F_{n-j}}B_{n,j+1}|)}{q_{n-j}}=0
\end{align}
Now $I_{F_{n-j}}I_{F_{n-j-1}}=I_{F_{n-j}},$ and hence
\begin{align}
|B_{n,j+1} - I_{F_{n-j}}B_{n,j+1}|=|(I_{F_{n-j-1}} -I_{F_{n-j}}) \prod_{\nu=j+1}^l \beta_{n,\nu}| \le \max \{1,c^l\}I_{F_{n-j-1} \cap F_{n-j}^c}.
\end{align}
and hence (2.7) implies (3.21).

Thus the induction holds, and to complete the proof of Lemma 3 it remains to verify (3.7). Hence let $\epsilon>0$ be given, and choose $k_{\epsilon} \ge 1$ and $c>0$ such that $\zeta(k) \le \epsilon$ for all $k \ge k_{\epsilon}$ and $\zeta(k) \le c$ for all $k \ge 1$. Then we have
\begin{align}
E(I_{ \{Z_n>0 \}} \zeta(Z_n)) \le \epsilon P(Z_n \ge k_{\epsilon}) +cP(1 \le Z_n \le k_{\epsilon}),
\end{align}
and (3.7) follows from (3.1) of Lemma 2. \QED

\bigskip

{\bf Proof of (3.2).} Set $\theta_0(B) = I_B(z)$ and $\theta_k(B)=P(H_k(\xi_1,\cdots,\xi_k) \in B)$ for $B$ a Borel subset of $M$ and $k \geq 1$, and let
$A_{n,i}=\{X_{n-i+1,Z_{n-i}} \in E_i\}$ for $ 1 \le i \le l<n.$ Then, 
$$
P(A_{n,i}|\mathcal{F}_{n-i})=\theta_{Z_{n-i}}(E_i) ~\rm {and}~ |E(I_{A_{n,i}} - \mu(E_i)|\mathcal{F}_{n-i})| = \zeta_i(Z_{n-i}),
$$
where $\zeta_i(k) = |\theta_k(E_i)-\mu(E_i)|$ for $ k \ge 0$ and $1 \le i \le l$. Since the sets $E_i$ are assumed to be continuity sets for the measure $\mu$, (2.6) implies
$$
\lim_{ k \rightarrow \infty} \zeta_i(k)=0,
$$
and hence (3.7) implies (3.3)  of Lemma 3 with $\gamma_{n,i}= I_{A_{n,i}}$, $\beta_{n,i}= \mu(E_i)$, and $1 \leq i \leq l$. Therefore, the assumptions of Lemma 3 hold, and hence (3.2) follows by taking $j=1$ in (3.6). Hence Theorem 1 holds.  \QED


\section{Proof of Theorem 2}

Application 1 of Theorem 1 implies (2.28) on $(C_0[0,1])^{\infty}$ with the product topology. Now we turn to its proof for the spaces $c_{0,\lambda}(C_0[0,1])$ and their stated norms $q_{\lambda}$. Given that weak convergence in the product topology implies the finite dimensional distributions of any finite set of coordinates converges in correct fashion, it suffices to show the probability measures of (2.27) are tight on the spaces $c_{0,\lambda}(C_0[0,1])$. This is the content of our next lemma. Its proof establishes Theorem 2.

\bigskip

\begin{lem}
Let $\{\mu_n:n \geq 1\}$ be as in (2.27), assume $m \geq 1$, and that $r(n) \rightarrow \infty$, with $r(n) = o(n)$. If (2.25) holds and $\lambda_j = (\delta_j \log(j+3))^{-\frac{1}{2}}$, where $\lim_{j \rightarrow \infty} \delta_j=\infty$, then the $\{\mu_n:n \geq 1\}$ are tight on $c_{0,\lambda}(C_0[0,1])$. Similarly, if (2.26) holds and $\lambda_j=j^{-\frac{(1+\delta)}{\rho}}$ for $\delta>0$, then we also have $\{\mu_n:n \geq 1\}$ tight on $c_{0,\lambda}(C_0[0,1])$.
\end{lem}

{\bf Proof}.  Since the finite dimensional distributions of any finite set of coordinates of $\{\mu_n\}$ converge weakly to the corresponding ones for $\mathcal{L}(B_1,B_2,\cdots)$ , standard arguments allow us to finish the proof by showing the $\{\mu_n\}$ are tight on $c_{0,\lambda}(C_0[0,1]))$.

To establish tightness we apply the remark in \cite{par}, p. 49. To show this remark applies we use the fact that the distributions of any finite set of coordinates are tight (since they are convergent), and therefore
it suffices to show for each $\epsilon >0$ that there exists a $d(\epsilon)$ such that $d \geq d(\epsilon)$ implies
\begin{align}
\limsup_{n \rightarrow \infty} P(q_{\lambda}(Q_d(\vxn)) \geq \epsilon|Z_{n-1}>0) \leq \epsilon.
\end{align}
Here $Q_d({\bf f}) =(0,\cdots,0,f_{d+1},f_{d+2},\cdots)$ for ${\bf f} \in (C_0[0,1])^{\infty}$. Since we are assuming $r(n)$ tends to infinity, for all n sufficiently large we have
$$
P(q_{\lambda}(Q_d(\vxn)) \geq \epsilon| Z_{n-1}>0) \leq \sum_{j=d+1}^{r(n)} I_{n,j},
$$
where
$$
I_{n,j}\equiv P(\max_{1 \leq l\leq Z_{n-j}} |\sum_{k=1}^{l} (\xi_{n-j+1,k} -m)| \geq Z_{n-j}^{\frac{1}{2}}\epsilon \lambda_j^{-1}|Z_{n-1} >0).
$$
Setting $J_{n,j} = I_{n,j}P(Z_{n-1}>0)$ we see
\begin{eqnarray}
J_{n,j} &=&\sum_{r=1}^{\infty}P(\max_{1 \leq l \leq r} |\sum_{k=1}^{l} (\xi_{n-j+1,k} -m)| \geq r^{\frac{1}{2}}\epsilon \lambda_j^{-1}, Z_{n-j}=r,Z_{n-1}>0).
\end{eqnarray}
Thus 
\begin{align}
J_{n,j}\le \sum_{r=1}^{\infty}P(\max_{1 \leq l \leq r} |\sum_{k=1}^{l} (\xi_{n-j+1,k} -m)| \geq r^{\frac{1}{2}}\epsilon \lambda_j^{-1}|Z_{n-j}=r)P(Z_{n-j}=r),
\end{align}
and by the branching property we see
\begin{align}
J_{n,j}\le \sum_{r=1}^{\infty}P(\max_{1 \leq l \leq r} |\sum_{k=1}^{l} (\xi_{k} -m)| \geq r^{\frac{1}{2}}\epsilon \lambda_j^{-1})P(Z_{n-j}=r),
\end{align}
where $\{\xi_k: k \ge 1\}$ are i.i.d. with law that of the offspring distribution. Since $\lambda_j^{-1} \rightarrow \infty$ there exists a $j_0=j_0(\epsilon)$ such that $j \geq j_0$ and Ottavianni's inequality implies
\begin{align}
J_{n,j}\le 2\sum_{r=1}^{\infty}P( |\sum_{k=1}^{r} (\xi_{k} -m)| \geq \frac{r^{\frac{1}{2}}\epsilon \lambda_j^{-1}}{2})P(Z_{n-j}=r).
\end{align}
Now under (2.25), Lemma 4.1 of \cite{KV} implies
for all $r \ge 1, j \ge 1$ that
$$
P( |\sum_{k=1}^{r} (\xi_{k} -m)| \geq \frac{r^{\frac{1}{2}}\epsilon \lambda_j^{-1}}{2})\leq 2 \exp\{-\theta \epsilon^2\lambda_j^{-2}/(64\beta)\}=2(j+3)^{-\frac{\theta \epsilon^2 \delta_j^2}{64\beta}}, 
$$
and hence for $j \geq j_0$ we have
\begin{align}
I_{n,j} \leq 4 (j+3)^{-\frac{\theta \epsilon^2 \delta_j^2}{64\beta}}P(Z_{n-j}>0)/P(Z_{n-1} >0).
\end{align}
When $m=1$ and $0< Var(Z_1)=\si^2< \infty$, we have by Theorem 1, p.19, of \cite{AN} that $\lim_{n \rightarrow \infty} nP(Z_n >0) = 2/ \si^2$. Hence for $n-j \geq n_0$ we have
\begin{align}
\frac{P(Z_{n-j} >0)}{P(Z_{n-1}>0)} = \frac{(n-j)P(Z_{n-j}>0)}{(n-1)P(Z_{n-1}>0)} \frac{(n-1)}{(n-j)} \leq2 \frac{(1 -\frac{1}{n})}{(1 - \frac{j}{n})}\leq \frac{2}{(1-\frac{j}{n})}.
\end{align}
Thus for $j=o(n), j \geq j_0,$ we have
\begin{align}
I_{n,j} \leq 16(j+3)^{-\frac{\theta \epsilon^2 \delta_j^2}{64\beta}}.
\end{align}
Now take $j_1= j_1(\epsilon)$ such that $j \geq j_1$ implies $\frac{\theta \epsilon^2 \delta_j^2}{64\beta}>2$. Given $\epsilon>0, r(n)= o(n)$ and $d>d_0(\epsilon) \equiv \max(j_0,j_1, \frac{16}{\epsilon} +1),$ we have
$$
\limsup_{n \rightarrow \infty}P(q_{\lambda}(Q_d(\vxn)) \geq \epsilon|Z_{n-1} >0) \leq \limsup_{n \rightarrow \infty} \sum_{j=d+1}^{r(n)} I_{n,j} \leq \epsilon. 
$$
Hence the lemma is proven under (2.25) if $m=1$. If $m>1$, then (4.8) is an even easier consequence of (4.6) since $\lim_{n \rightarrow \infty} P(Z_{n}>0) = 1-q>0$. Hence if $r(n)=o(n)$ , the lemma also holds in this case.

If (2.26) holds, then for all $r>0$ and $\rho\geq 2$ we have a constant $B_{\rho}<\infty$ such that an application of Markov's inequality and Corollary 8.2 in \cite{G}, p.151,
implies 
$$
P( |\sum_{k=1}^{r} (\xi_{k} -m)| \geq \frac{r^{\frac{1}{2}}\epsilon \lambda_j^{-1}}{2})\leq B_{\rho} \frac{(E(|\xi_1-m|^{\rho})}{(\frac{\epsilon \lambda_j^{-1}}{2})^{\rho}}.
$$
Hence the arguments can be completed as before, since under (2.26) we have $\lambda_j^{-1} = j^{\frac{(1+\delta)}{\rho}}$. Thus the lemma is proven. \QED

\section{Proof of Theorem 3}

Before we turn to the proof of Theorem 3 we provide a brief lemma, and recall that if $m>1$ and $~0<\si^2 \equiv E((Z_1-m)^2) < \infty$, then the Kesten-Stigum theorem, \cite{AN}, p. 9 (also see p. 24), implies that with probability one that
\begin{align}
\lim_{n \rightarrow \infty}W_n =W,
\end{align}
where $W_n= \frac{Z_n}{m^n}$, and
$W>0$ almost surely on the survival set $S$.


\begin{lem}
Under the given assumptions, we have almost surely that
\begin{align}
\lim_{n \rightarrow \infty} \frac{D_n}{m^n} = \sum_{j\geq 1} \frac{b_j}{m^j}W.
\end{align}
Furthermore, for $k \geq 1$ almost surely on S we have
\begin{align}
\lim_{n \rightarrow \infty} \frac{b_kZ_{n-k}}{D_n} = \theta_k,
\end{align}
where $\theta_k$ is given as in (2.29).
\end{lem}

{\bf Proof}. Observe that
$$
\frac{D_n}{m^n}=\sum_{k=1}^n \frac{b_ k Z_{n-k}}{m^n} = \sum_{k=1}^n \frac{b_k}{m^k} W + \sum_{k=1}^n \frac{b_k }{m^k}(\frac{Z_{n-k}}{m^{n-k}} - W),
$$
and since $\lim_{n\rightarrow \infty} Z_n/m^{n} =W $ almost everywhere by (5.1), an elementary argument combining $\kappa<\ff$ and $m>1$ easily implies  $\lim_{n\rightarrow \infty}  \sum_{k=1}^n \frac{b_k }{m^k}(\frac{Z_{n-k}}{m^{n-k}} - W)=0$ almost everywhere. Thus (5.2) holds. Combining (5.2) and (5.1) with $W>0$ almost surely on $S$, we thus have (5.3). Hence the lemma is proven. \QED

For the proof of Theorem 3 recall that if $D_n=0$, then we understand $X_n$ to be zero.
Furthermore,  if $D_n>0$, we then have 
$$
X_n = \sum_{j=1}^n \sqrt{\frac{b_j Z_{n-j}}{D_{n}}} \sqrt{b_j Z_{n-j}} (\frac{Z_{n-j+1}}{Z_{n-j}} -m),
$$
and for $1 \leq d \leq n$ we define
$$
X_{n,d} = \sum_{j=1}^d \sqrt{\frac{b_j Z_{n-j}}{D_{n}}} \sqrt{b_j Z_{n-j}} (\frac{Z_{n-j+1}}{Z_{n-j}} -m).
$$
Of course, when $D_n=0$, we understand $X_n$  and $X_{n,d}$ as given in these formulas to be zero.
We also use 
$$
\tilde X_n =\sqrt{\frac{D_n}{m^n}} X_n~ { \rm {and}} ~ \tilde X_{n,d} = \sqrt{\frac{D_n}{m^n}} X_{n,d},
$$
and their formulas analogous to those above for $X_n$ and $X_{n,d}$.

{\bf Proof of Theorem 3}. Take $\epsilon>0,$ and to simplify the notation set $\gamma_n=( \frac{D_n}{m^n})^{\frac{1}{2}}$. Then
\begin{align}
 P(X_n \leq x|S)=P(\tilde X_n \leq x \gamma_n|S),
 \end{align}
\begin{align}
P(\tilde X_n \leq x\gamma_n|S) \leq P(\tilde X_{n,d} \leq (x+\epsilon)\gamma_n|S) + P(|\tilde X_n - \tilde X_{n,d}| \geq \epsilon \gamma_n|S),
\end{align}
and
\begin{align}
P(|\tilde X_n - \tilde X_{n,d}| \geq \epsilon \gamma_n|S) \leq  P(|\tilde X_n - \tilde X_{n,d}| \geq \epsilon \delta|S) + P(0< \gamma_n< \delta|S).
\end{align}

Since $\epsilon>0$ is given, we choose $\delta >0$ sufficiently small that $P(0<W< 2 \delta^2)/P(S) < \epsilon$. Since $\lim_{n \rightarrow \infty} \gamma_n = W^{\frac{1}{2}}>0$ almost surely on $S$, there exists $n_0=n_0(\delta)$ such that $n\geq n_0$ implies
\begin{align}
P(0< \gamma_n< \delta|S) < \epsilon.
\end{align}

Once $\epsilon,\delta>0$ are fixed, we choose $d_0=d_0(\epsilon,\delta)$ such that $d \geq d_0$ implies that uniformly in $n$
\begin{align}
P(|\tilde X_n - \tilde X_{n,d}| \geq \epsilon \delta|S) \leq \epsilon.
\end{align}
To obtain $d_0$ we observe
that $P(|\tilde X_n - \tilde X_{n,d}| \geq \epsilon \delta|S) \leq P(|\tilde X_n - \tilde X_{n,d}| \geq \epsilon \delta)/P(S)$, and since $\kappa < \ff$ and the branching property easily implies
$$
E((\tilde X_n -\tilde X_{n,d})^2)= \sum_{j=d+1}^n m^{-n}b_jE(Z_{n-j-1})\si^2= \sum_{j=d+1}^n \si^2b_j m^{-(j+1)},
$$
we have $d_0=d_0(\epsilon,\delta)$ such that $d \ge d_0$ implies 
\begin{align}
E((\tilde X_n -\tilde X_{n,d})^2) < \epsilon^3 \delta^2.
\end{align}
Hence Markov's inequality, (5.9), and the above reasoning allows us to choose $d_0$ independent of $n$, so (5.8) holds.

Since $P( X_{n,d} \leq x+\epsilon|S)=P(\tilde X_{n,d} \leq (x+\epsilon)\gamma_n|S)$, by combining (5.4)-(5.8) we have for $d \geq d_0$ that
\begin{align}
P(X_n \leq x|S) \leq P( X_{n,d} \leq x+\epsilon|S) + 2\epsilon.
\end{align}
Similarly, we also have  for $d \geq d_0$ that
\begin{align}
P(X_n \leq x|S) \geq P( X_{n,d} \leq x-\epsilon|S) - 2\epsilon.
\end{align}
Now let
\begin{align}
X_{n,d}^{'} = \sum_{j=1}^d \sqrt{ \theta_j} \sqrt{b_j Z_{n-j}} (\frac{Z_{n-j+1}}{Z_{n-j}} -m),
\end{align}
and observe that by setting $t=1$ in the functional CLT of Application 1,  the continuous mapping theorem immediately implies
the uniform stochastic boundedness of 
$$
\{\sqrt{b_j Z_{n-j}} (\frac{Z_{n-j+1}}{Z_{n-j}} -m): 1 \leq j \leq d, n \geq 1 \}
$$
when these variables are conditioned on the event $\{Z_{n-1}>0\}$. Therefore, for each fixed $d$ we have from (5.3) of Lemma 5 and the previously mentioned uniform stochastic boundedness that 
\begin{align}
\lim_{n \rightarrow \infty} P(|X_{n,d} - X_{n,d}^{'}| \geq \epsilon|S)=0.
\end{align}
In addition, by the CLT  provided by Application 1 and the continuous mapping theorem we easily have
\begin{align}
\lim_{ n \rightarrow \infty} P(X_{n,d}^{'} \leq x|Z_{n-1}>0) = P(G_d \leq x)
\end{align}
for all real $x$, where $G_d$ is a mean zero Gaussian random variable with variance $\Lambda_d^2 =\sum_{j=1}^d \theta_j b_j\sigma^2$.

Since the events $\{Z_{n-1}>0\} \downarrow S$ with $P(S)=1-q>0$, by combining a standard argument implies (5.10), (5.13), and (5.14) we have for all $d \geq d_0$ and all real $x$ that
\begin{align}
\limsup_{n \rightarrow \infty}P(X_n \leq x|S) \leq P(G_d \leq x +2\epsilon) + 3\epsilon.
\end{align}
Using (5.11), a similar argument implies for all $d \geq d_0$ and all real $x$ that
\begin{align}
\liminf_{n \rightarrow \infty}P(X_n \leq x|S) \geq P(G_d \leq x -2\epsilon) - 3\epsilon.
\end{align}
Now take $d_1=d_1(\epsilon)$ sufficiently large such that $d\geq d_1$ implies  for all real $x$ that
\begin{align}
P(G_d \leq x) - \epsilon \leq P(G \leq x)  \leq P(G_d \leq x) +\epsilon,
\end{align}
where $G$ is as in the proposition. This condition follows easily since $\Lambda_d^2 \rightarrow \Lambda^2< \infty$. 

Letting $d$ tend to infinity in (5.15) and (5.16), (5,17) implies for all  $x$ that
\begin{align}
 \limsup_{n \rightarrow \infty}P(X_n \leq x|S) \leq P(G \leq x +2\epsilon) + 3\epsilon,
\end{align}
and
\begin{align}
 \liminf_{n \rightarrow \infty}P(X_n \leq x|S) \geq  P(G \leq x -2\epsilon) - 3\epsilon 
\end{align}
Letting $\epsilon \downarrow 0$ in (5.18) and (5.19), we have (2.30). Hence the the theorem is proven as the last claim is immediate from (2.30) . \QED


\section{Examples}

In this section we provide some examples where the CLT fails. We focus on the CLT as it is perhaps the result one might expect would be most likely to persist under suitable modifications of our basic assumptions. In the first example failure results from our branching process  $\{Z_n: n \geq 0\}$ being subcritical.  Hence, even though one has the same conditional independence structure as in the critical and supercritical cases, its 
behavior is quite different. In the other example the CLT fails through the use of deterministic normalizers. 
\bigskip

{\bf Subcritical Branching Fails the CLT}: Our result concerns the limit of 
$$
\mathcal{L}(Z_{n-1}^{\frac{1}{2}}(\frac{Z_n}{Z_{n-1}}-m)|Z_{n-1}>0),
$$
and shows that even for this
single distribution the CLT always fails. This is easy to see since the distribution of all the $\bar H_k \rm{'s}$ of the following lemma are discrete. 
\bigskip

\begin{lem} Assume that $E(Z_1^2) < \ff$ and set $L_n =Z_{n-1}^{\frac{1}{2}}(\frac{Z_n}{Z_{n-1}}-m)$. Then, for any $x \in\Real$,
\beq
\lim_{n \ra \ff} P(L_n \le x | Z_{n-1} >0) =  \sum_{k \ge 1} P(\sqrt{k}\bar{H}_k \le x) \theta_k,
\eeq
where $\bar{H}_k = \frac{1}{k}\sum_{i=1}^{k}(\xi_i-m) $,  $\{\xi_i:i \geq 1\}$ are i.i.d. with $\mathcal{L}(\xi_1)= \mathcal{L}(Z_1)$, and $\{\ta_k: k \geq 1\}$ is a probability distribution.
\end{lem}
\noindent{\bf Proof of Lemma 6}. Let $x \in \Real$. Then, by the branching property we easily have
\begin{eqnarray}
P(L_n \le x| Z_{n-1} >0) &=& \sum_{k \ge 0} P( L_n \le x ; Z_{n-1}=k | Z_{n-1} >0)\\
&=& \sum_{ k \ge 1} P( \sqrt{k}\bar{H}_k \le x) P(Z_{n-1}=k|Z_{n-1} >0).
\end{eqnarray}
Since $m<1$, Yaglom's Theorem on p.18 of \cite{AN} implies that $\lim_{n \rightarrow \infty} P(Z_{n-1} = k |Z_{n-1}>0)=\theta_k$, where $\{\theta_k: k \geq 1\}$ is a probability distribution. Thus, by the generalized dominated convergence theorem, it follows that
\begin{eqnarray}
\lim_{n \ra \ff}P(L_n \le x| Z_n >0) 
&=& \sum_{ k \ge 1} P( \sqrt{k}\bar{H}_k \le x) \theta_k.
\end{eqnarray}
This completes the proof of the lemma. \QED
\bigskip

{\bf Deterministic Normalizers Prevent the CLT}:
Even when $m>1$, the next example shows the spatial finite dimensional distributions related to Application 1 fail to be Gaussian when we use canonical deterministic normalizations $m^{\frac{n-1}{2}}$ instead of $Z_{n-1}^{\frac{1}{2}}$ in our CLT results. Of course, the motivation for these normalizations results from the Kesten-Stigum result, see (5.1), and in this situation the limit laws are a mixture of Gaussian laws and the random variable $W$ that appears in that result.

\begin{prop}
Let $m>1$, $0< \sigma^2 \equiv E((\xi_{1,1}-m)^2)< \infty$. For $i=1,\cdots,l<n$, set $H_{n-i+1}=\frac{m^{\frac{(n-i)}{2}}}{\sigma}(\frac{Z_{n-i+1}}{Z_{n-i}} -m)$ 
when $Z_{n-j}>0$, and zero otherwise, and let
\begin{align}
B_{n,i} = \{H_{n-i+1} \leq t_i \},
\end{align}
where $t_1,\cdots,t_l \in (-\infty,\infty)$ . Then
\begin{align}
\lim_{n \rightarrow \infty}  P(B_{n,1}\cap \cdots \cap B_{n,l}\cap \{Z_{n-1}>0\} )= 
E(\Phi(t_1 W^{1/2}) \cdots \Phi(t_l W^{1/2}) I_{S_0}),
\end{align}
where $S_0\equiv \{\lim_{n \rightarrow \infty} \frac{Z_n}{Z_{n-1}}=m\}$ and $W$ is as in (5.1).
\end{prop}
\bigskip


{\bf Proof of Proposition 1}. Let $\{\xi_j: j \ge 1\}$ be i.i.d. copies of $Z_1,$ and set
$$
h(0,t)=I_{[0,\infty)} ~\rm {and}~ h(k,t)=P(\sum_{i=1}^k(\xi_i-m) \le \sigma k^{\frac{1}{2}}t)
$$
for $k \ge 1$ and $t \in R^1$. Then $P(B_{n,i}|\mathcal {F}_{n-i}) =h(Z_{n-i},W_{n-i}t_i)$ and by the central limit theorem we have
$$
\zeta(k)\equiv \sup_{s \in R^1}|h(k,s)-\Phi(s)|  \rightarrow 0,
$$
where  $\Phi(\cdot)$ is the standard Gaussian cumulative distribution function. Therefore, we have
\begin{align}
|E(I_{B_{n,i}} - \Phi(W_{n-l}t_i)|\mathcal{F}_{n-i})| \le \zeta(Z_{n-i}) +|\Phi(W_{n-i}t_i)- \Phi(W_{n-l}t_i)|,
\end{align}
and since $m>1$, we have $P(Z_n >0) \downarrow 1-q=P(S_0)>0$ and (2.7) holding. Moreover, since $0<\si^2<\infty$ we have by (5.1) that $W_n \rightarrow W$ a.s. and $Z_n \rightarrow \infty$ on 
$S_0$. Hence the right hand side of (6.7) tends to zero on $S_0$ as $n$ tends to infinity, and the argument yielding (3.7), and $\{Z_{n>0}\} \downarrow S$ with $P(S\Delta S_0)=0$ implies for all $i=1,\cdots,l$ that
\begin{align}
E(|E(I_{B_{n,i}} - \Phi(W_{n-l}t_i)|\mathcal{F}_{n-i})|~ |Z_{n-i} >0) \rightarrow 0.
\end{align}
Hence we see that $\{\gamma_{n,i},\beta_{n,i}\} = \{I_{B_{n,i}},\Phi(W_{n-l}t_i)\} $ satisfies (3.3), as well as the remaining assumptions of Lemma 3. Therefore, Proposition 1 follows from (3.6) of Lemma 3. \QED


\end{document}